\documentclass[a4paper]{amsart}
\usepackage{amssymb}
\usepackage{verbatim}
\usepackage[dvips]{graphicx}
\input cyracc.def
\font\tencyr=wncyr10
\def\cyr{\tencyr\cyracc}
\theoremstyle{plain}
  \newtheorem{thm}{Theorem}[section]
  \newtheorem{lem}[thm]{Lemma}
  \newtheorem{cor}[thm]{Corollary}
  \newtheorem{prop}[thm]{Proposition}
  \newtheorem{conj}[thm]{Conjecture}

  \newtheorem*{obs*}{Observation}
\theoremstyle{definition}

\theoremstyle{remark}
  \newtheorem{rem}[thm]{Remark}
  \newtheorem*{ack}{Acknowledgments}
\newcommand{\Z}{\mathbb{Z}}
\newcommand{\C}{\mathbb{C}}
\newcommand{\R}{\mathbb{R}}
\newcommand{\Q}{\mathbb{Q}}

\newcommand{\Vol}{\operatorname{Vol}}
\newcommand{\CS}{\operatorname{CS}}

\newcommand{\Lob}{\operatorname{\text{\cyr L}}}
\newcommand{\Li}{\operatorname{Li}}

\newcommand{\length}{\operatorname{length}}
\newcommand{\torsion}{\operatorname{torsion}}

\renewcommand{\Re}{\operatorname{Re}}
\renewcommand{\Im}{\operatorname{Im}}
\numberwithin{equation}{section}
\allowdisplaybreaks
\begin{document}
\title[The colored Jones polynomials of the figure-eight knot]
{The colored Jones polynomials of the figure-eight knot
and its Dehn surgery spaces}
\author{Hitoshi Murakami}
\address{
Department of Mathematics,
Tokyo Institute of Technology,
Oh-okayama, Meguro, Tokyo 152-8551, Japan
}
\email{starshea@tky3.3web.ne.jp}
\author{Yoshiyuki Yokota}
\address{
Department of Mathematics,
Tokyo Metropolitan University,
Tokyo 192-0397, Japan}
\email{jojo@math.metro-u.ac.jp}
\date{\today}
\begin{abstract}
We calculate limits of the colored Jones polynomials of the figure-eight knot
and conclude that in most cases they determine the volumes and the Chern--Simons
invariants of the three-manifolds obtained by Dehn surgeries along it.
\end{abstract}
\thanks{This research is partially supported by Grant-in-Aid for Scientific
Research (C) (15540089) and (B) (15340019).}
\maketitle
\section{Introduction}
Let $K$ be a knot and $J_N(K;t)$ be the colored Jones polynomial of $K$ in
the three-sphere $S^3$ corresponding to the $N$-dimensional representation of
$sl_2(\C)$ (see for example \cite{Kirby/Melvin:INVEM91}).
We normalize it so that $J_N(\text{unknot};t)=1$.
\par
In \cite{Kashaev:LETMP97}, R.~Kashaev introduced a series of numerical link
invariants and proposed a conjecture that a limit of his invariants would
determines the hyperbolic volume of the complement of a knot if it has a
complete hyperbolic structure.
On the other hand in \cite{Murakami/Murakami:ACTAM101} J.~Murakami and the first
author proved that Kashaev's invariant is equal to the absolute value of
$J_N\bigl(K;\exp(2\pi\sqrt{-1}/N)\bigr)$.
Moreover in \cite{Murakami/Murakami/Okamoto/Takata/Yokota:EXPMA02} J.~Murakami,
M.~Okamoto, T.~Takata and the authors studies several knots and proposed the
following complexification of Kashaev's conjecture.
\begin{conj}
Let $K$ be a hyperbolic knot (that is, the complement of $K$ possesses a
complete hyperbolic structure).
Then
\begin{equation*}
  2\pi\lim_{N\to\infty}\frac{\log{J_N\bigl(K;\exp(2\pi\sqrt{-1}/N)\bigr)}}{N}
  =
  \Vol(S^3\setminus{K})+\sqrt{-1}\CS(S^3\setminus{K}),
\end{equation*}
where $\Vol$ is the hyperbolic volume and $\CS$ is the Chern--Simons
invariant for knots \cite{Meyerhoff:density}.
\end{conj}
See
\cite{Murakami:ALDT_VI,
Murakami:Nagoya99,
D.Thurston:Grenoble,
Yokota:Murasugi70,
Yokota:volume00,
Yokota:Topology_Symposium2000,
Yokota:GTM02,
Yokota:INTIS03}
for related topics.
\par
On the other hand the first author studied in \cite{Murakami:KYUMJ2004} other limits
of the colored Jones polynomials of the figure-eight knot and showed the
following theorem.
\begin{thm}\label{thm:real}
Let $E$ be the figure-eight knot and $r$ a irrational number satisfying
${5/6}<{r}<{7/6}$ or $r=1$.
Then
\begin{equation*}
  2\pi{r}\lim_{N\to\infty}
  \frac{\log\left|J_N\left(E;\exp(2\pi{r}\sqrt{-1}/N)\right)\right|}{N}
\end{equation*}
is the volume of the cone manifold with cone angle $2\pi|1-r|$
whose singularity is the figure-eight knot.
\end{thm}
\par
He also calculated
\cite{Murakami:SURIK00}
`fake' limits (optimistic limits) of the Witten--Reshetikhin--Turaev
invariants of three-manifolds obtained from $S^3$ by integral
Dehn surgeries along the figure-eight knot and observed they give the volumes
and the Chern--Simons invariants of those manifolds.
\par
The purpose of this paper is to show that limits of the colored
Jones polynomials of the figure-eignt give the volume and the Chern--Simons
invariants of the three-manifold obtained by general surgeries.
To do that we study the asymptotic behavior of $J_N(E;\exp(2\pi{r}\sqrt{-1}/N))$
for large $N$ by using the saddle point method with $E$ the figure-eight knot.
Then it turns out that it determines the potential function introduced by
W.~Neumann and D.~Zagier \cite{Neumann/Zagier:TOPOL85} to study the volumes of
Dehn surgered spaces along knots.
(See also \cite{Gukov:COMMP2005}.)
The we apply T.~Yoshida's theorem \cite{Yoshida:INVEM85} to show that the limit
of $\log\Bigl(J_N\bigl(E;\exp(2\pi{r}\sqrt{-1}/N)\bigl)\Bigr)/N$
($N\to\infty$) gives the volume and the Chern--Simons invariant when $r\in\C$ is
close to $1$.
\begin{ack} Part of this work was done during the authors were attending the workshops `Physics and Geometry of three-dimensional Quantum Gravity' at the International Centre for Mathematical Sciences (ICMS) in Edinburgh, 29 June-5 July 2003, and `Hyperbolic volume conjecture' at University of Geneva, 8-14 September 2003. The authors thank the organizers and the participants. \par Thanks are also due to Hiroshi Isozaki and Rinat Kashaev for their useful conversations. \end{ack} 
\section{Main result}
To state our main result we recall some of the results in
\cite{Neumann/Zagier:TOPOL85,Yoshida:INVEM85}.
\par
Let $K$ be a hyperbolic knot in the three-sphere $S^3$, that is, the complement
$S^3\setminus{K}$ has a complete hyperbolic structure.
The complete hyperbolic structure can be deformed by a complex parameter
$u$ that is the logarithm of the ratio of the eigenvalues of the image of the meridian $K$ by the holonomy representation
$\pi_1(S^3\setminus{K})\to\mathrm{SL}_2(\C)$.
We also denote by $v$ the logarithm of the eigenvalue of the image of the
longitude.
If the structure is complete, then $u=0$.
Let $\Phi(u)$ be the potential function introduced in
\cite[Theorem 3]{Neumann/Zagier:TOPOL85} that satisfies
\begin{equation*}
\begin{cases}
  \dfrac{d\,\Phi(u)}{du}=2v,
  \\[3mm]
  \Phi(0)=0.
\end{cases}
\end{equation*}
We also put
\begin{equation*}
  f(u):=\frac{\Phi(u)}{4}-\frac{uv}{4}.
\end{equation*}
(see \cite[Equation (50)]{Neumann/Zagier:TOPOL85}).
\par
Let $\ell$ be a loop on the boundary of the tubular neighborhood
$N(K)$ of $K$ in $S^3$.
Suppose that the hyperbolic structure parameterized by $u$ on $S^3\setminus{K}$
can be completed to an orbifold $K_{u}$ by attaching a solid
torus such that its meridian coincides with $\ell$.
It is called generalized Dehn surgery and if $pu+qv=2\pi\sqrt{-1}$ for
coprime integers $p$ and $q$, then it is usual $(p,q)$-Dehn surgery and
$K_{u}$ becomes a closed three-manifold
\cite[Chapter 4]{Thurston:GT3M} and $\ell$ is homotopic to
$p\cdot\text{meridian}+q\cdot\text{longitude}$ in $\pi_1(\partial{N(K)})$.
Let $\gamma$ be the core of the solid torus.
It is a geodesic in $K_{u}$, and its length and torsion are denoted by
$\length(\gamma)$ and $\torsion(\gamma)$ respectively.
(If one travels around $\gamma$ then the plane perpendicular to it is
rotated by $\torsion(\gamma)$.
See \cite[Definition~1.1]{Yoshida:INVEM85} for the precise definition.)
Then the volume and the Chern--Simons invariant of $K_{u}$ can be obtained
as follows
(\cite[Theorem 2]{Yoshida:INVEM85}, which was conjectured in
\cite[Conjecture]{Neumann/Zagier:TOPOL85}).
\begin{thm}[T.~Yoshida]\label{thm:Yoshida}
Put $\lambda(\gamma):=\length(\gamma)+\sqrt{-1}\torsion(\gamma)$.
Then we have
\begin{multline*}
  \Vol(K_{u})+\sqrt{-1}\CS(K_{u})
  \equiv
  \frac{f(u)}{\sqrt{-1}}
  -\frac{\pi}{2}\lambda(\gamma)
  +\Vol\left(S^3\setminus{K}\right)+\sqrt{-1}\CS\left(S^3\setminus{K}\right)
  \\
  \mod{\pi^2\sqrt{-1}\Z}.
\end{multline*}
\end{thm}
\par
Now we study the case where the knot is the figure-eight knot $E$.
\par
Let $u$ be a complex parameter near $0$.
Put
\begin{align*}
  m&:=-\exp(u/2),
  \\
  z&:=\frac{-m^4+m^2+1+\sqrt{(m^2+m+1)(m^2+m-1)(m^2-m+1)(m^2-m-1)}}{2m^2},
  \\
  \intertext{and}
  w&:=\frac{m^4+m^2-1+\sqrt{(m^2+m+1)(m^2+m-1)(m^2-m+1)(m^2-m-1)}}{2m^2}.
\end{align*}
We take the branch of the square root so that $z=w=\exp(\pi/3)$ when
$u=0$.
We also put
\begin{equation*}
  v:=2\log\bigl(z(1-z)\bigr),
\end{equation*}
where the branch is chosen so that $v=0$ when $u=0$.
Note that $z$ and $w$ satisfy the following equations.
\begin{equation*}
  \begin{cases}
    \log{w}+\log\bigl(1-z\bigr)=u,
    \\
    \log{z}+\log\bigl(1-z\bigr)+\log{w}+\log\bigl(1-w\bigr)=0.
  \end{cases}
\end{equation*}
\par
Put
\begin{align*}
  f(u)&:=\frac{1}{2\pi}\left\{R(z)+R(w)-\frac{\pi^2}{6}\right\}
        -\frac{\sqrt{-1}}{2\pi}\Vol\left(S^3\setminus{E}\right)
  \\
  \intertext{and}
  \Phi(u)&:=4f(u)+uv,
\end{align*}
where
\begin{equation*}
  R(\xi):=\frac{1}{2}\log{\xi}\log(1-\xi)
  -\int_{0}^{\xi}\frac{\log(1-\eta)}{\eta}d\eta
\end{equation*}
is Roger's dilogarithm function.
Then
\begin{gather*}
  \begin{cases}
    \dfrac{d\,\Phi(u)}{du}=2v,
    \\[3mm]
    \Phi(0)=0,
  \end{cases}
  \\
  \intertext{and}
  \Vol(E_u)+\sqrt{-1}\CS(E_u)
  \equiv
  \frac{f(u)}{\sqrt{-1}}-\frac{\pi}{2}\lambda(\gamma)
  +\Vol\left(S^3\setminus{E}\right)
  \mod{\pi^2\sqrt{-1}\Z},
\end{gather*}
since $\CS\left(S^3\setminus{E}\right)=0$.
\par
Now we state our main result.
\par
Let $J_N(K;t)$ be the colored Jones polynomial of a knot $K$ corresponding
to the $N$-dimensional representation of $sl_2(\C)$.
We normalize it so that $J_N(\text{unknot};t)=1$.
Then we have the following theorem.
\begin{thm}\label{thm:main}
Let $E$ be the figure-eight knot.
There exists a neighborhood $U$ of $0$ in $\C$ such that for any
$u\in({U}\setminus\pi\sqrt{-1}\Q)\cup\{0\}$, the following limit exists:
\begin{equation*}
  (u+2\pi\sqrt{-1})
  \lim_{N\to\infty}
  \frac{\log\Bigl(J_N\bigl(E;\exp\left((u+2\pi\sqrt{-1})/N\right)\bigr)\Bigr)}
  {N}.
\end{equation*}
Moreover if we denote the limit by $H(u)$, then it can be extended to a complex
analytic function and it satisfies the following equalities.
\begin{equation*}
\begin{cases}
  \dfrac{d\,H(u)}{du}=\dfrac{v}{2}+\pi\sqrt{-1},
  \\[3mm]
  H(0)=\sqrt{-1}\Vol(S^3\setminus{E}).
\end{cases}
\end{equation*}
\end{thm}
\begin{rem}
The precise meaning of the limit is as follows.
If the ratio
\begin{equation*}
  \frac{J_N\bigl(E;\exp\left((u+2\pi\sqrt{-1})/N\right)\bigr)}
  {\exp\bigl(N\,h(u)\bigr)}
\end{equation*}
grows polynomially with respect to $N$, then we denote the limit
\begin{equation*}
  \frac{\log\Bigl(J_N\bigl(E;\exp\left((u+2\pi\sqrt{-1})/N\right)\bigr)\Bigr)}
  {N}
\end{equation*}
by $h(u)$.
So it is defined modulo $2\pi\sqrt{-1}$.
\end{rem}
As a corollary we can express the volume and the Chern--Simons invariant
of a three-manifold obtained from the figure-eight knot by generalized Dehn
surgery.
\begin{cor}\label{cor:Vol_CS}
Let $E_u$ be the three-manifold obtained from the figure-eight knot
by generalized Dehn surgery corresponding to $u$ near $0$ such that
$\pi{u}\notin\Q$.
Then
\begin{equation*}
  \Vol(E_u)+\sqrt{-1}\CS(E_u)
  \equiv
  \frac{H(u)}{\sqrt{-1}}-\pi{u}-\frac{uv}{4\sqrt{-1}}
  -\frac{\pi}{2}\lambda(\gamma)
  \mod{\pi^2\sqrt{-1}\Z}.
\end{equation*}
\end{cor}
\section{Limits of the colored Jones polynomials}
K.~Habiro \cite{Habiro:SURIK00} and T.~Le showed that the colored Jones
polynomial $J_N(E;t)$ of the figure-eight knot $E$ is given as follows
(see also \cite{Masbaum:AGT03}).
\begin{equation}\label{eq:Jones}
  J_N(E;t)=\sum_{n=0}^{N-1}\prod_{k=1}^{n}
           t^N\left(1-t^{-N-k}\right)\left(1-t^{-N+k}\right).
\end{equation}
We will study the asymptotic behavior of
$J_N(E;q_r)$ with $q_r:=\exp(2\pi{r}\sqrt{-1}/N)$ for large $N$.
\par
We approximate $\prod_{k=1}^{n}\left(1-q_r^{-N{\pm}k}\right)$ by
an integral.
Note that
\begin{equation*}
  \prod_{k=1}^{n}\left(1-q_r^{-N{\pm}k}\right)
  =
  \exp
  \left\{
    \frac{N}{2\pi}
    \sum_{k=1}^{n}
    \frac{2\pi}{N}
    \log\bigl(1-\exp(\pm2\pi{k}{r}\sqrt{-1}/N-2\pi{r}\sqrt{-1})\bigr)
  \right\}.
\end{equation*}
\begin{lem}\label{lem:max_min}
If $r\not\in\R$ then
\begin{gather*}
  \log\left|1-e^{2\pi{b}/N}\right|
  \le
  \log\left|1-\exp\left(sr\sqrt{-1}-2\pi{r}\sqrt{-1}\right)\right|
  \le
  \log\left(1+e^{2\pi{b}}\right),
  \\
  \intertext{and}
  \log\left|1-e^{2\pi{b}}\right|
  \le
  \log\left|1-\exp\left(-sr\sqrt{-1}-2\pi{r}\sqrt{-1}\right)\right|
  <
  \log\left(1+e^{4\pi|b|}\right)
\end{gather*}
where $r=a+b\sqrt{-1}$ with $a\in\R$, $0\ne{b}\in\R$ and
$0\le{s}\le2\pi-2\pi/N$.
\begin{proof}
If we put $r:=a+b\sqrt{-1}$ with $b\ne0$, we have
\begin{equation*}
  \exp\left(sr\sqrt{-1}-2\pi{r}\sqrt{-1}\right)
  =
  e^{-bs+2\pi{b}}\exp\bigl((as-2\pi{a})\sqrt{-1}\bigr).
\end{equation*}
Therefore
$\left|1-\exp\left(sr\sqrt{-1}-2\pi{r}\sqrt{-1}\right)\right|$
is smaller (greater, respectively) than or equal to the longest
(shortest, respectively) distance between $1\in\C$
and the circle centered at $0$ with radius $e^{-bs+2\pi{b}}$.
So if $b>0$,
\begin{equation*}
  1+e^{2\pi{b}}
  \ge
  \left|1-\exp\left(sr\sqrt{-1}-2\pi{r}\sqrt{-1}\right)\right|
  \ge
  e^{-b(2\pi-2\pi/N)+2\pi{b}}-1
  =
  e^{2\pi{b}/N}-1
\end{equation*}
and if $b<0$,
\begin{equation*}
  1+e^{2\pi{b}}
  \ge
  \left|1-\exp\left(sr\sqrt{-1}-2\pi{r}\sqrt{-1}\right)\right|
  \ge
  1-e^{2\pi{b}/N}.
\end{equation*}
\par
Similarly if $b>0$,
\begin{equation*}
  e^{2\pi{b}}-1
  \le
  \left|1-\exp\left(-sr\sqrt{-1}-2\pi{r}\sqrt{-1}\right)\right|
  \le
  e^{b(2\pi-2\pi/N)+2\pi{b}}+1
  <
  e^{4\pi{b}}+1
\end{equation*}
and if $b<0$,
\begin{equation*}
  1-e^{2\pi{b}}
  \le
  \left|1-\exp\left(-sr\sqrt{-1}-2\pi{r}\sqrt{-1}\right)\right|
  \le
  e^{4\pi{b}-2\pi{b}/N}+1
  \le
  e^{2\pi{b}}+1.
\end{equation*}
\end{proof}
\end{lem}
We put
\begin{multline*}
  \varphi_{N,\pm}(n)
  :=
  \sum_{k=1}^{n}
  \frac{2\pi}{N}
  \log
  \bigl(
    1-\exp(\pm2\pi{k}r\sqrt{-1}/N-2\pi{r}\sqrt{-1})
  \bigr)
  \\
  -
  \int_{0}^{2{\pi}n/N}
  \log
  \bigl(
    1-\exp({\pm}sr\sqrt{-1}-2{\pi}r\sqrt{-1})
  \bigr)\,ds.
\end{multline*}
\begin{lem}\label{lem:logN}
If $r\not\in\R$, then there exist constants $d\ge1$ and $\delta>0$ such that
\begin{gather*}
  \left|\Re\varphi_{N,+}(n)\right|
  <
  \frac{2\pi{d}}{N}
  \left\{
    \log\left(1+e^{2\pi{b}}\right)
    -
    \log\left(\frac{2\pi|b|\delta}{N}\right)
  \right\}
  \\
  \intertext{and}
  \left|\Re\varphi_{N,-}(n)\right|
  <
  \frac{2\pi{d}}{N}
  \left\{
    \log\left(1+e^{4\pi|b|}\right)
    -
    \log\left|1-e^{2\pi{b}}\right|
  \right\}
\end{gather*}
for sufficiently large $N$.
\end{lem}
\begin{proof}
First we note that if a real continuous function $g(s)$ ($0\le{s}\le{A}$) has $c$ extremes, then
\begin{equation*}
  \left|
    \int_{0}^{nA/N}g(s)\,ds
    -
    \sum_{k=1}^{n}\frac{A}{N}g(kA/N)
  \right|
  \le
  (2c+1)\frac{A}{N}
  \left(\max_{0\le{s}\le{A}}{g(s)}-\min_{0\le{s}\le{A}}{g(s)}\right).
\end{equation*}
This is because if $g$ is monotonic in $[jA/N,(j+l)A/N]$
\begin{multline*}
  \left|
    \int_{jA/N}^{(j+l)A/N}g(s)\,ds
    -
    \sum_{k=j+1}^{j+l}\frac{A}{N}g(kA/N)
  \right|
  \\
  \le
  \frac{A}{N}
  \left|
    g(jA/N)-g((j+l)A/N)
  \right|
  \le
  \frac{A}{N}
  \left(
    \max_{0\le{s}\le{A}}{g(s)}-\min_{0\le{s}\le{A}}{g(s)}
  \right)
\end{multline*}
and there are $c+1$ such intervals.
For other $c$ intervals $[jA/N,(j+1)A/N]$ clearly
\begin{equation*}
  \left|
    \int_{jA/N}^{(j+1)A/N}g(s)\,ds
    -
    \frac{A}{N}g((j+1)A/N)
  \right|
  \le
  \frac{A}{N}
  \left(
    \max_{0\le{s}\le{A}}{g(s)}-\min_{0\le{s}\le{A}}{g(s)}
  \right)
\end{equation*}
holds.
\par
Therefore from Lemma~\ref{lem:max_min} there exists a constant $d$ such that
\begin{gather*}
  \left|
    \Re\varphi_{N,+}(n)
  \right|
  \le
  \frac{2\pi{d}}{N}
  \left\{
    \log\left(1+e^{2\pi{b}}\right)
    -
    \log\left|1-e^{2\pi{b}/N}\right|
  \right\}
  \\
  \intertext{and}
  \left|\Re\varphi_{N,-}(n)\right|
  <
  \frac{2\pi{d}}{N}
  \left\{
    \log\left(1+e^{4\pi|b|}\right)
    -
    \log\left|1-e^{2\pi{b}}\right|
  \right\}
\end{gather*}
Since if $b>0$,
\begin{equation*}
  e^{2\pi{b}/N}-1
  =
  \frac{2\pi{b}}{N}+\frac{(2\pi{b})^2}{2N^2}+\cdots
  >
  \frac{2\pi{b}}{N}
\end{equation*}
for any $N$ and if $b<0$,
\begin{equation*}
  1-e^{2\pi{b}/N}
  =
  -\frac{2\pi{b}}{N}-\frac{(2\pi{b})^2}{2N^2}-\cdots
  >
  -\frac{2\pi{b}\delta}{N}
\end{equation*}
for any $0<\delta<1$ if $N$ is sufficiently large.
This completes the proof.
\end{proof}
Therefore we have
\begin{equation*}
  D^{-1}N^{-d}
  <
  \left|
    \exp\left(\frac{N}{2\pi}\varphi_{N,\pm}(n)\right)
  \right|
  <
  DN^d
\end{equation*}
for some constants $D>0$ and $d\ge1$ for sufficient large $N$.
\par
We also have
\begin{lem}
If $r\not\in\R$, then for any positive number $\varepsilon$
\begin{equation*}
  \left|\Im\varphi_{N,\pm}(n)\right|
  <
  \varepsilon
\end{equation*}
for sufficiently large $N$.
\end{lem}
\begin{proof}
Clearly
\begin{multline*}
  \max_{2\pi/N\le{s}\le2\pi-2\pi/N}
  \arg\bigl(1-\exp(\pm{sr}\sqrt{-1}-2\pi{r}\sqrt{-1})\bigr)
  \\
  -
  \min_{2\pi/N\le{s}\le2\pi-2\pi/N}
  \arg\bigl(1-\exp(\pm{sr}\sqrt{-1}-2\pi{r}\sqrt{-1})\bigr)
\end{multline*}
is bounded by a positive constant depending only $r$.
So we can show as in the previous lemma that
\begin{multline*}
  \left|\Im\varphi_{N,\pm}(n)\right|
  =
  \Biggl|
    \sum_{k=1}^{n}\frac{2\pi}{N}
    \arg\bigl(1-\exp(\pm{sr}\sqrt{-1}-2\pi{r}\sqrt{-1})\bigr)
  \\
    -
    \int_{0}^{2\pi{n}/N}
    \arg\bigl(1-\exp(\pm{sr}\sqrt{-1}-2\pi{r}\sqrt{-1})\bigr)
  \Biggr|
\end{multline*}
is bounded by a positive constant times $1/N$ and the result follows.
\end{proof}
Putting
\begin{equation*}
  \chi_{N,\pm}(n):=\exp\left(\frac{N}{2\pi}\varphi_{N,\pm}(n)\right),
\end{equation*}
we have
\begin{multline*}
  \prod_{k=1}^{n}
  \left(1-q_r^{-N\pm{k}}\right)
  \\
  =
  \chi_{N,\pm}(n)
  \exp
  \left\{
    \frac{N}{2\pi}
    \int_{0}^{2\pi{n}/N}
    \log\bigl(1-\exp(\pm{sr}\sqrt{-1}-2\pi{r}\sqrt{-1})\bigr)\,ds
  \right\}
\end{multline*}
with $D^{-1}N^{-d}<\left|\chi_{N,\pm}(n)\right|<DN^d$ and
$|\arg\chi_{N,\pm(n)}|<\varepsilon$ if $r\not\in\R$.
Therefore
\begin{equation*}
\begin{split}
  &\prod_{k=1}^{n}
  \left(1-q_r^{-N\pm{k}}\right)
  \\
  &\quad=
  \chi_{N,\pm}(n)
  \exp
  \left\{
    \frac{\pm{N}}{2\pi{r}\sqrt{-1}}
    \int_{\exp(-2\pi{r}\sqrt{-1})}
        ^{\exp(\pm2\pi{nr}\sqrt{-1}/N-2\pi{r}\sqrt{-1})}
    \frac{\log(1-u)}{u}\,du
  \right\}
  \\
  &\quad=
  \chi_{N,\pm}(n)
  \exp
  \left\{
    \frac{\pm{N}}{2\pi{r}\sqrt{-1}}
    \Bigl(
      \Li_2\left(m^{-2}\right)-\Li_2\left(q_r^{\pm{n}}m^{-2}\right)
    \Bigr)
  \right\}.
\end{split}
\end{equation*}
Here we put
\begin{equation*}
  m:=\exp(\pi{r}\sqrt{-1})
\end{equation*}
and
\begin{equation*}
  \Li_2(z):=-\int_{0}^{z}\frac{\log(1-u)}{u}\,du
\end{equation*}
is the dilogarithm function.
\par
Thus we have
\begin{equation*}
\begin{split}
  &J_{N}(E;q_r)
  \\
  &\quad=
  \sum_{n=0}^{N-1}
  q_r^{Nn}
  \prod_{k=1}^{n}
  \left(1-q_r^{-N+k}\right)
  \left(1-q_r^{-N-k}\right)
  \\
  &\quad=
  \sum_{n=0}^{N-1}
  q_r^{Nn}
  \chi_{N,+}(n)
  \exp
  \left\{
    \frac{N}{2\pi{r}\sqrt{-1}}
    \Bigl(
      \Li_2\left(m^{-2}\right)-\Li_2\left(q_r^{n}m^{-2}\right)
    \Bigr)
  \right\}
  \\
  &\qquad\qquad\qquad\times
  \chi_{N,-}(n)
  \exp
  \left\{
    \frac{-N}{2\pi{r}\sqrt{-1}}
    \Bigl(
      \Li_2\left(m^{-2}\right)-\Li_2\left(q_r^{-n}m^{-2}\right)
    \Bigr)
  \right\}
  \\
  &\quad=
  \sum_{n=0}^{N-1}
  \chi_N(n)
  \exp
  \left\{
    \frac{N}{2\pi{r}\sqrt{-1}}
    \Bigl(
      \Li_2\left(q_r^{-n}m^{-2}\right)-\Li_2\left(q_r^{n}m^{-2}\right)
  \right.
  \\
  &\qquad\qquad\qquad\qquad
  \left.\vphantom{\frac{N}{2\pi{r}\sqrt{-1}}}
      +\bigl(\log\left(-q_r^n\right)+\pi\sqrt{-1}\bigr)(\log{m^2}
    \Bigr)
  \right\}
  \\
  &\quad=
  \sum_{n=0}^{N-1}
  \chi_N(n)
  \exp
  \left\{
    \frac{N}{2\pi{r}\sqrt{-1}}H\left(q_r^n,m^2\right)
  \right\},
\end{split}
\end{equation*}
where we put
\begin{equation*}
  H(\xi,\eta)
  :=
  \Li_2\left(\xi^{-1}\eta^{-1}\right)-\Li_2\left(\xi\eta^{-1}\right)
  +\left(\log(-\xi)+\pi\sqrt{-1}\right)\log{\eta}
\end{equation*}
and $\chi_N(n):=\chi_{N,+}(n)\chi_{N,-}(n)$.
Note that $D^{-2}N^{-2d}<|\chi(n)|<D^2N^{2d}$.
Note also that we use $\log(-\xi)+\pi\sqrt{-1}$ instead of $\log\xi$ since we choose the branch of $\log$ as $(0,+\infty)$ and that of $\Li_2$ as $(1,+\infty)$.
\par
Let $\Psi_N(z)$ be an analytic function such that $\Psi_N\bigl(\exp(2\pi{n}\sqrt{-1}/N)\bigr)=\chi_N(n)$, that $P(N^{-1})<|\Psi_N(z)|<Q(N)$ with some polynomials $P$ and $Q$, and that $|\arg\Psi_N(z)|<\varepsilon$ near $C$ defined below. Since the residue around $z=\exp(2\pi{n}\sqrt{-1}/N)$ of $1/\bigl(z(1-z^{-N})\bigr)$ is $1/N$, we have from the residue theorem
\begin{equation*}   J_N(E;q_r)   =   1+   \frac{N}{2\pi\sqrt{-1}}   \int_{C}   \frac{\Psi_N(z)}{z(1-z^{-N})}   \exp   \left(\frac{N}{2\pi{r}\sqrt{-1}}H(z^r,m^2)\right)   dz, \end{equation*} where $C=C_{+}{\cup}C_{-}$ is defined by
\begin{equation*}
\begin{cases}
  C_{+}&:=
  \left\{
    z\in\C\mid|z|=1+\kappa,\pi/N\le\arg{z}\le2\pi-\pi/N
  \right\}
  \\
  &\phantom{:=}\quad
  \cup
  \left\{
    t\exp( \sqrt{-1}/N)\in\C\mid1-\kappa\le{t}\le1+\kappa
  \right\}
  \\
  &\phantom{:=}\quad
  \cup
  \left\{
    t\exp(-\sqrt{-1}/N)\in\C\mid1-\kappa\le{t}\le1+\kappa
  \right\}
  \\
  C_{-}&:=
  \left\{
    z\in\C\mid|z|=1-\kappa,\pi/N\le\arg{z}\le2\pi-\pi/N
  \right\}
\end{cases}
\end{equation*}
for small $\kappa>0$.
Now we have
\begin{equation*}
\begin{split}
  &J_N(E;q_r)
  \\   &\quad=
  1+
  \frac{N}{2\pi\sqrt{-1}}
  \int_{C_{+}}
  \frac{\Psi_N(z)}{z}
  \exp
  \left(\frac{N}{2\pi{r}\sqrt{-1}}H(z^r,m^2)\right)
  dz
  \\
  &\quad\quad
  +
  \frac{N}{2\pi\sqrt{-1}}
  \int_{C_{+}}
  \frac{\Psi_N(z)}{z(1-z^{-N})}
  \exp
  \left(
    \frac{N}{2\pi{r}\sqrt{-1}}
    \left\{
      H(z^r,m^2)
      -
      2\pi\sqrt{-1}\log{z^r}
    \right\}
  \right)
  dz
  \\
  &\quad\quad
  +
  \frac{N}{2\pi\sqrt{-1}}
  \int_{C_{-}}
  \frac{\Psi_N(z)}{z(z^{N}-1)}
  \exp
  \left(
    \frac{N}{2\pi{r}\sqrt{-1}}
    \left\{
      H(z^r,m^2)
      +
      2\pi\sqrt{-1}\log{z^r}
    \right\}
  \right)
  dz.
\end{split}
\end{equation*}
\par
Here we want to apply the steepest descent method
(or the saddle point method, see for example
\cite[Theorem 7.2.9]{Marsden/Hoffman:Complex_Analysis})
\begin{thm}[Steepest Descent Method]\label{thm:steepest_descent_method}
Let $\Gamma\colon[a,b]\to\C$ be a $C^1$ curve.
Let $h(z)$ be a continuous function on $\Gamma$ that is analytic at
$z_0=\Gamma(t_0)$, and $g_N(z)$ be a continuous function on $z\in{\Gamma}$
for each $N=1,2,3,\dots$ with $g_N(z_0)\ne0$ such that $|g_N(z)|<N^p$ for some
$p>0$.
We also assume the following conditions:
For sufficient large integer $N$,
\begin{enumerate}
\item
  $\int_{\Gamma}g_N(z)\exp\bigl(Nh(z)\bigr)dz$ converges absolutely.
\item
  $d\,h(z_0)/dz=0$ and $d^2\,h(z_0)/dz^2\ne0$.
\item
  $\Im\bigl(Nh(z)\bigr)$ is constant for $z$ on $c$ in some neighborhood of
  $z_0$.
\item
  $\Re\bigl(Nh(z)\bigr)$ takes its strict maximum along $\Gamma$ at $z_0$.
\end{enumerate}
Then
\begin{equation*}
  \int_{\Gamma}g_N(z)\exp\bigl(Nh(z)\bigr)dz
  \underset{N\to\infty}{\sim}
  \frac{\sqrt{2\pi}g_N(z_0)\exp\bigl(Nh(z_0)\bigr)}
       {\sqrt{N}\sqrt{-d^2\,h(z_0)/dz^2}}
\end{equation*}
for appropriate chosen square roots.
Here $\underset{N\to\infty}{\sim}$ means the ratio of both hand sides
goes to $1$ when $N\to\infty$.
\end{thm}
\begin{rem}
In 
\cite[Theorem 7.2.9]{Marsden/Hoffman:Complex_Analysis}
it is assumed that $g_N(z)$ does not depend on $N$ but the result is
also true if we assume that $|g_N(z)|$ grows polynomially.
See Exercises 8 and 10 in Section 7.3 of
\cite{{Marsden/Hoffman:Complex_Analysis}}.
\end{rem}
Putting $h(z):=\frac{1}{2\pi{r}\sqrt{-1}}H(z^r,m^2)$,
$g(z):=\Psi_N(z)/z$, we can prove the following proposition.
\begin{prop}
We have
\begin{multline*}
  \int_{C_{+}}
  \frac{\Psi_N(z)}{z}
  \exp
  \left(\frac{N}{2\pi{r}\sqrt{-1}}H(z^r,m^2)\right)
  dz
  \\
  \underset{N\to\infty}{\sim}
  \frac{\sqrt{2\pi}\Psi_N(y^{1/r})}{r\sqrt{N}\sqrt{y-y^{-1}}}
  \exp\left(\frac{N}{2\pi{r}\sqrt{-1}}H(y,m^2)\right)
\end{multline*}
with
\begin{equation*}
  y:=\frac{m^4-m^2+1-\sqrt{(m^2+m+1)(m^2+m-1)(m^2-m+1)(m^2-m-1)}}{2m^2},
\end{equation*}
where the branch of the square root in the definition of $y$ is taken
so that if $m=-1$, then $y=\exp(-\pi\sqrt{-1}/3)$.
Note that $y$ satisfies the following equality:
\begin{equation}\label{eq:y}
  y+y^{-1}=m^2-1+m^{-2}.
\end{equation}
\end{prop}
\begin{proof}
First of all we calculate derivatives of $H(z^r,m^2)$.
We have
\begin{equation*}
\begin{split}
  \frac{\partial\,H(z^r,m^2)}{\partial{z}}
  &=
  \frac{\partial\,H(\xi,\eta)}{\partial\xi}
  \frac{d\,z^r}{dz}
  \\
  &=
  \frac{rz^{r-1}}{z^r}
  \left\{
    \log(1-z^{-r}m^{-2})
    +
    \log(1-z^{r}m^{-2})
    +
    \log{m^2}
  \right\}
  \\
  &=
  \frac{r}{z}
  \log\left(m^2-(z^{r}+z^{-r})+m^{-2}\right)
\end{split}
\end{equation*}
and
\begin{multline*}
  \frac{\partial^2\,H(z^r,m^2)}{\partial{z^2}}
  \\
  =
  -
  \frac{r}{z^2}\log\left(m^2-(z^{r}+z^{-r})+m^{-2}\right)
  +
  \frac{r}{z}
  \times
  \frac{-rz^{r-1}+rz^{-r-1}}{m^2-(z^{r}+z^{-r}+m^{-2})}.
\end{multline*}
Therefore $y^{1/r}$ gives a solution to $\partial\,H(z^r,m^2)/\partial{z}=0$,
and
\begin{equation*}
  \frac{\partial^2\,H(z^r,m^2)}{\partial{z^2}}\Biggr|_{z=y^{1/r}}
  =
  r^2y^{-2/r}(-y+y^{-1}).
\end{equation*}
\par
Now we will consider the case when $r=1$ precisely.
In this case $y=\exp(-\pi\sqrt{-1}/3)$.
We change $C_{+}$ to ${C'}_{+}$ where
\begin{equation*}
\begin{split}
  {C'}_{+}&:=
  \left\{
    z\in\C\mid|z|=1,\pi/N\le\arg{z}\le2\pi-\pi/N
  \right\}
  \\
  &\phantom{:=}\quad
  \cup
  \left\{
    t\exp( \sqrt{-1}/N)\in\C\mid1-\kappa\le{t}\le1
  \right\}
  \\
  &\phantom{:=}\quad
  \cup
  \left\{
    t\exp(-\sqrt{-1}/N)\in\C\mid1-\kappa\le{t}\le1
  \right\}.
\end{split}
\end{equation*}
Then on the unit circle
\begin{equation}\label{eq:H_on_circle}
\begin{split}
  H(\exp(\theta\sqrt{-1}),1)
  &=
  \Li_2(\exp(-\theta\sqrt{-1}))-\Li_2(\exp(\theta\sqrt{-1}))
  \\
  &=
  \left\{
    \frac{\pi^2}{6}
    -
    \frac{2\pi-\theta}{2}
    \left(\pi-\frac{2\pi-\theta}{2}\right)
    +
    2\sqrt{-1}\Lob\left(\frac{2\pi-\theta}{2}\right)
  \right\}
  \\
  &\quad
  -
  \left\{
    \frac{\pi^2}{6}
    -
    \frac{\theta}{2}\left(\pi-\frac{\theta}{2}\right)
    +
    2\sqrt{-1}\Lob\left(\frac{\theta}{2}\right)
  \right\}
  \\
  &=
  -4\sqrt{-1}\Lob\left(\frac{\theta}{2}\right).
\end{split}
\end{equation}
Here we use the following formula:
\begin{equation}\label{eq:dilog}
  \Li_2\bigl(\exp(\beta\sqrt{-1})\bigr)
  =
  \frac{\pi^2}{6}
  -
  \frac{\beta}{2}\left(\pi-\frac{\beta}{2}\right)
  +
  2\sqrt{-1}\Lob\left(\frac{\beta}{2}\right)
\end{equation}
for $0\le\beta\le2\pi$ and the fact that $\Lob$ is an odd function with period
$\pi$.
(See for example \cite[p. 18]{Milnor:BULAM382}.)
So on ${C'}_{+}$ $\Im{H(z,1)}$ takes its maximum at $y$ and near $y$
$\Re{H(z,1)}=0$.
Moreover in this case
\begin{equation*}
  \frac{\partial^2\,H(z,1)}{\partial{z^2}}\Biggr|_{z=y}
  =
  \frac{-3-\sqrt{3}\sqrt{-1}}{2}
  \ne0.
\end{equation*}
Therefore we can apply the saddle point method and the proposition follows when
$r=1$.
\par
If $r$ is near $1$ we can also apply the saddle point method similarly from the
continuity of $H(z^r,m^2)$ with respect to $r$, completing the proof.
\end{proof}
\par
From Lemma~\ref{lem:C_-} below we see
\begin{equation*}
\begin{split}
  &
  \left|
  \int_{C_{+}}
  \frac{\Psi_N(z)}{z(1-z^{-N})}
  \exp
  \left(
    \frac{N}{2\pi{r}\sqrt{-1}}
    \left\{
      H(z^r,m^2)
      -
      2\pi\sqrt{-1}\log{z^r}
    \right\}
  \right)
  dz
  \right|
  \\
  &\quad
  \le
  \int_{C_{+}}
  \frac{\left|\Psi_N(z)\right|}{\left|z(1-z^{-N})\right|}
  \exp
  \left(
    \frac{N}{2\pi}
    \Im
    \left(
      \frac{H(z^r,m^2)-2\pi\sqrt{-1}\log{z^r}}{r}
    \right)
  \right)
  |dz|
  \\
  &\quad
  <
  \exp
  \left(
    \frac{N}{2\pi}
    \Im\left(\frac{H(y,m^2)}{r}\right)
  \right)
  \int_{C_{+}}
  \frac{\left|\Psi_N(z)\right|}{\left|z(1-z^{-N})\right|}
  |dz|
  \\
  &\quad
  \le
  \exp
  \left(
    \frac{N}{2\pi}
    \Im\left(\frac{H(y,m^2)}{r}\right)
  \right)
  \frac{|C_{+}|Q(N)}{(1+\kappa)\left(1-(1+\kappa)^{-N}\right)}
\end{split}
\end{equation*}
and
\begin{multline*}
  \left|
    \int_{C_{-}}
    \frac{\Psi_N(z)}{z(z^{N}-1)}
    \exp
    \left(
      \frac{N}{2\pi{r}\sqrt{-1}}
      \left\{
        H(z^r,m^2)
        +
        2\pi\sqrt{-1}\log{z^r}
      \right\}
    \right)
    dz
  \right|
  \\
  <
  \exp
  \left(
    \frac{N}{2\pi}
    \Im\left(\frac{H(y,m^2)}{r}\right)
  \right)
  \frac{|C_{-}|Q(N)}{(1-\kappa)\left(1-(1-\kappa)^{N}\right)},
\end{multline*}
where $|C_{\pm}|$ is the length of $C_{\pm}$.
Now we have
\begin{equation}\label{eq:log/N}
\begin{split}
  &
  \log\bigl(J_N(E;q_r)\bigr)\big/N
  -
  \frac{H(y,m^2)}{2\pi{r}\sqrt{-1}}
  -
  \log\left(
        \frac{\sqrt{2\pi}\Psi_N(y^{1/r})}{r\sqrt{N}\sqrt{y-y^{-1}}}
      \right)\Big/N
  \\
  &\quad
  =
  \log
  \Biggl(
    \frac{r\sqrt{N}\sqrt{y-y^{-1}}}{\sqrt{2\pi}\Psi_N(y^{1/r})}
    \exp\left(\frac{-N}{2\pi{r}\sqrt{-1}}H(y,m^2)\right)
  \\
  &\quad\quad
  +
  \frac{N}{2\pi\sqrt{-1}}
  \frac{\int_{C_+}\frac{\Psi_N(z)}{z}
        \exp\left(\frac{N}{2\pi{r}\sqrt{-1}}H(z^r,m^2)\right)\dz}
       {\frac{\sqrt{2\pi}\Psi_N(y^{1/r})}{r\sqrt{N}\sqrt{y-y^{-1}}}
         \exp\left(
               \frac{N}{2\pi{r}\sqrt{-1}}H(y,m^2)
             \right)}
  \\
  &\quad\quad
  +
  \frac{N}{2\pi\sqrt{-1}}
  \frac{\int_{C_+}\frac{\Psi_N(z)}{z(1-z^{-N})}
        \exp
        \left(
          \frac{N}{2\pi{r}\sqrt{-1}}H(z^r,m^2)-2\pi\sqrt{-1}\log{z^r}
        \right)\dz}
       {\frac{\sqrt{2\pi}\Psi_N(y^{1/r})}{r\sqrt{N}\sqrt{y-y^{-1}}}
         \exp\left(
               \frac{N}{2\pi{r}\sqrt{-1}}H(y,m^2)
             \right)}
  \\
  &\quad\quad
  +
    \frac{N}{2\pi\sqrt{-1}}
    \frac{\int_{C_-}\frac{\Psi_N(z)}{z(z^{N}-1)}
          \exp
          \left(
            \frac{N}{2\pi{r}\sqrt{-1}}H(z^r,m^2)+2\pi\sqrt{-1}\log{z^r}
          \right)\dz}
         {\frac{\sqrt{2\pi}\Psi_N(y^{1/r})}{r\sqrt{N}\sqrt{y-y^{-1}}}
           \exp\left(
                 \frac{N}{2\pi{r}\sqrt{-1}}H(y,m^2)
               \right)}
  \Biggr)
  \Bigg/N.
\end{split}
\end{equation}
We will estimate each term in $\log$.
First we observe that
\begin{multline*}
  \left|
    \frac{r\sqrt{N}\sqrt{y-y^{-1}}}{\sqrt{2\pi}\Psi_N(y^{1/r})}
    \exp\left(\frac{-N}{2\pi{r}\sqrt{-1}}H(y,m^2)\right)
  \right|
  \\
  =
  \frac{|r|\sqrt{N}|\sqrt{y-y^{-1}}|}{\sqrt{2\pi}\left|\Psi_N(y^{1/r})\right|}
  \exp\left(-N\Re\left(\frac{H(y,m^2)}{2\pi{r}\sqrt{-1}}\right)\right)
\end{multline*}
goes to $0$.
This is because $P(N^{-1})<|\Psi_N(y^{1/r})|<Q(N)$ for polynomials $P$ and $Q$,
and $\Re\bigl(H(y,m^2)/2\pi{r}\sqrt{-1}\bigr)>0$ if $r$ is near $1$ since
when $r=1$
\begin{equation*}
  \frac{H(y,m^2)}{2\pi{r}\sqrt{-1}}
  =
  \frac{H(\exp(-\pi\sqrt{-1}/3),1)}{2\pi\sqrt{-1}}
  =
  \frac{2\Lob(\pi/6)}{\pi}
  >
  0
\end{equation*}
from \eqref{eq:H_on_circle}.
Next we have
\begin{equation*}
  \lim_{N\to\infty}
  \frac{\int_{C_+}\frac{\Psi_N(z)}{z}
        \exp\left(\frac{N}{2\pi{r}\sqrt{-1}}H(z^r,m^2)\right)\dz}
         {\frac{\sqrt{2\pi}\Psi_N(y^{1/r})}{r\sqrt{N}\sqrt{y-y^{-1}}}
           \exp\left(
                 \frac{N}{2\pi{r}\sqrt{-1}}H(y,m^2)
               \right)}
  =1,
\end{equation*}
\begin{equation*}
\begin{split}
  &
  \left|
    \frac{N}{2\pi\sqrt{-1}}
    \frac{\int_{C_+}\frac{\Psi_N(z)}{z(1-z^{-N})}
          \exp
          \left(
            \frac{N}{2\pi{r}\sqrt{-1}}H(z^r,m^2)-2\pi\sqrt{-1}\log{z^r}
          \right)\dz}
         {\frac{\sqrt{2\pi}\Psi_N(y^{1/r})}{r\sqrt{N}\sqrt{y-y^{-1}}}
           \exp\left(
                 \frac{N}{2\pi{r}\sqrt{-1}}H(y,m^2)
               \right)}
  \right|
  \\
  &\quad
  <
  \frac{|rN^{3/2}\sqrt{y-y^{-1}}||C_{+}|Q(N)}
        {(2\pi)^{3/2}|\Psi_N(y^{1/r}|)
         (1+\kappa)\left(1-(1+\kappa)^{-N}\right)}
  \\
  &\quad
  <
  \frac{|rN^{3/2}\sqrt{y-y^{-1}}||C_{+}|Q(N)(1+N\kappa)}
        {(2\pi)^{3/2}P(N^{-1}))(1+\kappa)N\kappa}
\end{split}
\end{equation*}
and
\begin{equation*}
\begin{split}
  &
  \left|
    \frac{N}{2\pi\sqrt{-1}}
    \frac{\int_{C_-}\frac{\Psi_N(z)}{z(z^{N}-1)}
          \exp
          \left(
            \frac{N}{2\pi{r}\sqrt{-1}}H(z^r,m^2)+2\pi\sqrt{-1}\log{z^r}
          \right)\dz}
         {\frac{\sqrt{2\pi}\Psi_N(y^{1/r})}{r\sqrt{N}\sqrt{y-y^{-1}}}
           \exp\left(
                 \frac{N}{2\pi{r}\sqrt{-1}}H(y,m^2)
               \right)}
  \right|
  \\
  &\quad
  <
  \frac{|rN^{3/2}\sqrt{y-y^{-1}}||C_{-}|Q(N)}
       {(2\pi)^{3/2}|\Psi_N(y^{1/r}|)
         (1-\kappa)\left(1-(1-\kappa)^{N}\right)}
  \\
  &\quad
  <
  \frac{|rN^{3/2}\sqrt{y-y^{-1}}||C_{-}|Q(N)}
       {(2\pi)^{3/2}P(N^{-2})(1-\kappa)(\kappa+(N-1)\kappa^2/2)}.
\end{split}
\end{equation*}
So we see that the absolute values of the first, third, and fourth terms in
$\log$ in \eqref{eq:log/N} are smaller than a polynomial in $N$.
Therefore \eqref{eq:log/N} goes to $0$ when $N\to\infty$ and so we have
\begin{equation*}
\begin{split}
  &\lim_{N\to\infty}\frac{\log\bigl(J_N(E;q_r)\bigr)}{N}
  \\
  &\quad=
  \frac{H(y,m^2)}{2\pi{r}\sqrt{-1}}
  +
  \lim_{N\to\infty}
  \log
  \left(\frac{\sqrt{2\pi}\Psi_N(y^{1/r})}{r\sqrt{N}\sqrt{y-y^{-1}}}\right)
  \big/N
  \\
  &\quad=
  \frac{H(y,m^2)}{2\pi{r}\sqrt{-1}}
  +
  \lim_{N\to\infty}\frac{\log\left|\Psi_N(y^{1/r})\right|}{N}
  +
  \sqrt{-1}\lim_{N\to\infty}\frac{\arg\Psi_N(y^{1/r})}{N}
  \\
  &\quad=
  \frac{H(y,m^2)}{2\pi{r}\sqrt{-1}}
\end{split}
\end{equation*}
since $P(N^{-1})<|\Psi_N(y^{1/r})|<Q(N)$ and
$|\arg\Psi_N(y^{1/r})|<\varepsilon$.
\begin{lem}\label{lem:C_-}
There exists a neighborhood $U'$ of $1\in\C$ such that if $r{\in}U'$,
\begin{equation}\label{eq:ineq}
  \Im
  \left(
    \frac{H\left(z^r,m^2\right)\mp2\pi\sqrt{-1}\log{z^r}}{r}
  \right)
  <
  \Im
  \left(
    \frac{H(y,m^2)}{r}
  \right)
\end{equation}
for $z\in{C_{\pm}}$.
\end{lem}
\begin{proof}
We will show that the inequality holds when $r=1$.
Then from the compactness of $C_{\pm}$ and the continuity we see that there
exists such a neighborhood $U'$ of $1$.
\par
Put
\begin{equation*}
\begin{split}
  h_{\mp}(z)
  &:=
  H(z,1)\mp2\pi\sqrt{-1}\log{z}-H\bigl(\exp(-\pi\sqrt{-1}/3),1\bigr)
  \\
  &=
  \Li_2(1/z)-\Li_2(z)\mp2\pi\sqrt{-1}\log{z}
  -\Li_2\left(e^{\pi\sqrt{-1}/3}\right)
  +\Li_2\left(e^{-\pi\sqrt{-1}/3}\right).
\end{split}
\end{equation*}
and we regard it as a function of $\rho$ and $\theta$ with
$z=\rho\exp(\theta\sqrt{-1})$; that is
\begin{multline*}
  h_{\mp}(\rho,\theta)
  =
  \Li_2\left(\rho^{-1}{e}^{-\theta\sqrt{-1}}\right)
  -
  \Li_2\left(\rho{e}^{\theta\sqrt{-1}}\right)
  \pm
  2\pi\theta
  \mp
  2\pi\sqrt{-1}\log{\rho}
  \\
  -\Li_2\left(e^{\pi\sqrt{-1}/3}\right)
  +\Li_2\left(e^{-\pi\sqrt{-1}/3}\right).
\end{multline*}
We also put
\begin{equation*}
\begin{cases}
  D_{+}:=\{z\in\C\mid|z|\ge1\},
  \\
  D_{-}:=\{z\in\C\mid|z|\le1\}.
\end{cases}
\end{equation*}
We will show that $\Im{h_{\mp}(\rho,\theta)}$ on $D_{\pm}$ takes its unique
maximum $0$ at $\exp(-\pi\sqrt{-1}/3)$.
\par
Since
\begin{equation*}
  \frac{\partial\,h_{\mp}(\rho,\theta)}{\partial\theta}
  =
  \sqrt{-1}\log\left(1-\rho^{-1}{e}^{-\theta\sqrt{-1}}\right)
  +
  \sqrt{-1}\log\left(1-\rho{e}^{\theta\sqrt{-1}}\right)
  \pm
  2\pi,
\end{equation*}
we have
\begin{equation*}
  \frac{\partial\,\Im{h_{\mp}(\rho,\theta)}}{\partial\,\theta}
  =
  \log
  \biggl|
    \left(1-\rho^{-1}{e}^{-\theta\sqrt{-1}}\right)
    \left(1-\rho     {e}^{ \theta\sqrt{-1}}\right)
  \biggr|.
\end{equation*}
Since
\begin{equation*}
\begin{split}
  \left|
    \left(1-\rho^{-1}{e}^{-\theta\sqrt{-1}}\right)
    \left(1-\rho     {e}^{ \theta\sqrt{-1}}\right)
  \right|^2
  &=
  \left(1+\rho^{-2}-2\rho^{-1}\cos\theta\right)
  \left(1+\rho^{ 2}-2\rho     \cos\theta\right)
  \\
  &=
  \left(\rho+\rho^{-1}-2\cos\theta\right)^2,
\end{split}
\end{equation*}
$\Im{h_{\mp}(\rho,\theta)}$ takes its maximum at
$\theta=2\pi-\arccos\left(\dfrac{\rho+\rho^{-1}-1}{2}\right)$ for a fixed $\rho$
if $\left(3-\sqrt{5}\right)/2<\rho<\left(3+\sqrt{5}\right)/2$.
\par
Now we put $\theta(\rho):=2\pi-\arccos\left(\dfrac{\rho+\rho^{-1}-1}{2}\right)$
and denote $h_{\mp}(\rho,\theta(\rho))$ by $h_{\mp}(\rho)$.
Since
\begin{equation*}
  \frac{\partial\,h_{\mp}(\rho,\theta)}{d\rho}
  =
  \frac{1}{\rho}
  \left\{
    \log\left(1-\rho^{-1}{e}^{-\theta\sqrt{-1}}\right)
    +
    \log\left(1-\rho{e}^{\theta\sqrt{-1}}\right)
    \mp
    2\pi\sqrt{-1}
  \right\}
\end{equation*}
and
\begin{equation*}
  \frac{d\,\theta(\rho)}{d\rho}
  =
  \frac{1-\rho^{-2}}{2\sin{\theta(\rho)}}
\end{equation*}
we have
\begin{multline*}
  \frac{d\,h_{\mp}(\rho)}{d\rho}
  =
  \left\{
    \frac{1}{\rho}+\frac{\sqrt{-1}(1-\rho^{-2})}{2\sin\theta(\rho)}
  \right\}
  \\
  \times
  \left\{
    \log
    \biggl(
      \left(1-\rho^{-1}{e}^{-\theta(\rho)\sqrt{-1}}\right)
      \left(1-\rho     {e}^{ \theta(\rho)\sqrt{-1}}\right)
    \biggr)
    \mp2\pi\sqrt{-1}
  \right\}.
\end{multline*}
Therefore
\begin{equation*}
\begin{split}
  \frac{d\,\Im{h_{\mp}(\rho)}}{d\rho}
  &=
  \frac{1}{\rho}
  \left\{
    \arg
    \biggl(
      \left(1-\rho^{-1}{e}^{-\theta(\rho)\sqrt{-1}}\right)
      \left(1-\rho     {e}^{ \theta(\rho)\sqrt{-1}}\right)
    \biggr)
    \mp2\pi
  \right\}
  \\
  &\quad
  +
  \frac{1-\rho^{-2}}{2\sin\theta(\rho)}
  \left\{
    \log
    \biggl|
      \left(1-\rho^{-1}{e}^{-\theta(\rho)\sqrt{-1}}\right)
      \left(1-\rho     {e}^{ \theta(\rho)\sqrt{-1}}\right)
    \biggr|
  \right\}
  \\
  &=
  \frac{1}{\rho}
  \left\{
    \arg
    \biggl(
      \left(1-\rho^{-1}{e}^{-\theta(\rho)\sqrt{-1}}\right)
      \left(1-\rho     {e}^{ \theta(\rho)\sqrt{-1}}\right)
    \biggr)
    \mp2\pi
  \right\}
\end{split}
\end{equation*}
since
\begin{equation*}
  \biggl|
    \left(1-\rho^{-1}{e}^{-\theta(\rho)\sqrt{-1}}\right)
    \left(1-\rho     {e}^{ \theta(\rho)\sqrt{-1}}\right)
  \biggr|
  =1.
\end{equation*}
\par
Thus we have
\begin{equation*}
\begin{cases}
  \dfrac{d\,\Im{h_{+}(\rho)}}{d\rho}>0,
  \\[3mm]
  \dfrac{d\,\Im{h_{-}(\rho)}}{d\rho}<0
\end{cases}
\end{equation*}
and so $\Im{h_{\pm}(\rho)}$ on $D_{\mp}$ takes maximum at $\rho=1$.
\par
Since $\theta(1)=-\pi/3$, $\Im{h_{\mp}(1)}=0$.
Therefore \eqref{eq:ineq} holds when $r=1$.
\end{proof}
\par
Now we consider the case when $r\in\R$.
In this case $m^2+m^{-2}=2\cos(2\pi{r})$ and so
\begin{equation*}
  y=\exp\bigl(-\alpha(r)\sqrt{-1}\bigr)
\end{equation*}
with $\alpha(r):=\arccos(\cos(2\pi{r})-1/2)$ ($0\le\arccos{x}\le{\pi}$).
From \eqref{eq:dilog}
\begin{equation*}
\begin{split}
  &
  H(y,m^2)
  \\
  &\quad=
  \Li_2(y^{-1}m^{-2})-\Li_2(ym^{-2})+\left(\log(-y)+\pi\sqrt{-1}\right)\log{m^2}
  \\
  &\quad=
  \Li_2
  \Bigl(\exp\bigl(\alpha(r)-2\pi{r})\sqrt{-1}\bigr)\Bigr)
  -
  \Li_2
  \Bigl(\exp\bigl(-\alpha(r)-2\pi{r})\sqrt{-1}\bigr)\Bigr)
  \\
  &\quad\quad
  +2\pi(r-1)(\alpha(r)-2\pi)
  \\
  &\quad
  =
  \frac{\pi^2}{6}
  -\frac{2\pi+\alpha(r)-2\pi{r}}{2}
   \left(\pi-\frac{2\pi+\alpha(r)-2\pi{r}}{2}\right)
  +2\sqrt{-1}\Lob\left(\frac{\alpha(r)-2\pi{r}}{2}\right)
  \\
  &\quad\quad
  -\frac{\pi^2}{6}
  +\frac{4\pi-\alpha(r)-2\pi{r}}{2}
   \left(\pi-\frac{4\pi-\alpha(r)-2\pi{r}}{2}\right)
  -2\sqrt{-1}\Lob\left(\frac{-\alpha(r)-2\pi{r}}{2}\right)
  \\
  &\quad\quad
  +2\pi(r-1)(\alpha(r)-2\pi)
  \\
  &\quad=
  -2\pi^2(r-1)
  +
  2\sqrt{-1}
  \left\{
    \Lob\bigl(\pi{r}+\alpha(r)/2\bigr)-\Lob\bigl(\pi{r}-\alpha(r)/2\bigr)
  \right\}.
\end{split}
\end{equation*}
\par
On the other hand from \cite[Theorem 1.2]{Murakami:KYUMJ2004}
\begin{multline*}
  2\pi{r}\Re\lim_{N\to\infty}
  \frac{\log{J_N\bigl(E;\exp(2\pi{r}\sqrt{-1}/N)\bigr)}}{N}
  \\
  =
  2
  \left\{
    \Lob\bigl(\pi{r}+\alpha(r)/2\bigr)-\Lob\bigl(\pi{r}-\alpha(r)/2\bigr)
  \right\}
\end{multline*}
if $r\not\in\Q$ and $5/6<r<7/6$, or $r=1$.
Moreover from its proof we see that the sign of
$J_N\bigl(E;\exp(2\pi{r}\sqrt{-1})/N\bigr)$ for large $N$ is
$(-1)^{\lfloor{N(1-r)/r}\rfloor}$, where $\lfloor{x}\rfloor$ is the greatest
integer that does not exceed $x$.
So
\begin{equation*}
\begin{split}
  2\pi{r}\Im\lim_{N\to\infty}
  \frac{\log{J_N\bigl(E;\exp(2\pi{r}\sqrt{-1}/N)\bigr)}}{N}
  &=
  2\pi{r}\Im\lim_{N\to\infty}
  \frac{\log\left((-1)^{\lfloor{N(1-r)/r}\rfloor}\right)}{N}
  \\
  &=
  2\pi^2(1-r).
\end{split}
\end{equation*}
Therefore we have
\begin{equation*}
  2\pi{r}\sqrt{-1}\lim_{N\to\infty}
  \frac{\log{J_N\bigl(E;\exp(2\pi{r}\sqrt{-1}/N)\bigr)}}{N}
  =
  H(y,m^2)
\end{equation*}
if $5/6<r<7/6$ and $r\not\in\Q\setminus\{0\}$.
\par
Putting
\begin{equation*}
  u:=2\pi{r}\sqrt{-1}-2\pi\sqrt{-1}
\end{equation*}
and
\begin{equation*}
  H(u):=H(y,m^2),
\end{equation*}
we have proved the following theorem.
\begin{thm}
Let $E$ be the figure-eight knot.
There exists a neighborhood $U$ of $0$ in $\C$ such that for any
$u\in(U\setminus\pi\sqrt{-1}\Q)\cup\{0\}$,
\begin{equation*}
  (u+2\pi\sqrt{-1})
  \lim_{N\to\infty}
  \frac{\log\left(
              J_N\left(
                   E;\exp\left(\frac{u+2\pi\sqrt{-1}}{N}\right)
                 \right)
            \right)}
  {N}
  =H(u).
\end{equation*}
\end{thm}
\section{Proof of the main theorem}
In this section we will prove the main theorem.
\par
First we will calculate the derivative of $H(u):=H(y,m^2)$.
\begin{lem}
\begin{equation*}
  \frac{d\,\bigl(H(u)-\pi\sqrt{-1}u\bigr)}{du}
  =
  \frac{v}{2}.
\end{equation*}
\end{lem}
\begin{proof}
We first calculate the partial derivatives of $H(\xi,\eta)$.
We have
\begin{equation*}
\begin{split}
  \frac{\partial\,H(\xi,\eta)}{\partial\xi}
  &=
  \
  \frac{1}{\xi}
  \left\{
     \log\left(1-\xi^{-1}\eta^{-1}\right)
    +\log\left(1-\xi\eta^{-1}\right)
    +\log\eta
  \right\}
  \\
  &=
  \frac{1}{\xi}\log\left(\eta+\eta^{-1}-\xi-\xi^{-1}\right)
\end{split}
\end{equation*}
and
\begin{equation*}
\begin{split}
  \frac{\partial\,H(\xi,\eta)}{\partial\eta}
  &=
  \frac{1}{\eta}
  \left\{
     \log\left(1-\xi^{-1}\eta^{-1}\right)
    -\log\left(1-\xi\eta^{-1}\right)
    +\log(-\xi)
  \right\}
  \\
  &=
  \frac{1}{\eta}
  \left\{
    \log\left(\frac{1-\xi\eta}{\eta-\xi}\right)
    +\pi\sqrt{-1}
  \right\}.
\end{split}
\end{equation*}
Since
\begin{equation*}
  \left.\frac{\partial\,H(\xi,\eta)}{\partial\,\xi}\right|_{\xi=y,\eta=m^2}
  =
  \frac{1}{y}\log\left(m+m^{-1}-y-y^{-1}\right)
  =0
\end{equation*}
from \eqref{eq:y} we have
\begin{equation*}
\begin{split}
  \frac{d\,H(u)}{du}
  &=
  \left.\frac{\partial\,H(\xi,\eta)}{\partial\xi}\right|_{\xi=y,\eta=m^2}
  \times\frac{d\,y}{du}
  +
  \left.\frac{\partial\,H(\xi,\eta)}{\partial\eta}\right|_{\xi=y,\eta=m^2}
  \times\frac{d\,m^2}{du}
  \\
  &=
  \log\left(\frac{1-ym^2}{m^2-y}\right)+\pi\sqrt{-1}.
\end{split}
\end{equation*}
But since
\begin{equation*}
\begin{cases}
  ym^2+zm^2=1
  \\
  z^{-1}+m^2y^{-1}=1,
\end{cases}
\end{equation*}
we have
\begin{equation*}
  \frac{d\,\bigl(H(u)-\pi\sqrt{-1}u\bigr)}{du}
  =
  \log\bigl(z(1-z)\bigr)
  =
  \frac{v}{2}.
\end{equation*}
\end{proof}
\par
Next we will calculate $H(0)$.
\begin{lem}
\begin{equation*}
  H(0)=\sqrt{-1}\Vol\left(S^3\setminus{E}\right)
\end{equation*}
\end{lem}
\begin{proof}
Since $m\bigr|_{u=0}=-1$ and $y\bigr|_{u=0}=\exp(-\pi\sqrt{-1}/3)$, we have
\begin{equation*}
  H(0)
  =
   \Li_2\left(e^{-\pi\sqrt{-1}/3}\right)
  -\Li_2\left(e^{ \pi\sqrt{-1}/3}\right)
  =
  4\sqrt{-1}\Lob(\pi/6)
\end{equation*}
from \eqref{eq:H_on_circle}.
\par
On the other hand from the equality
\begin{equation*}
  \Lob(2\theta)=2\Lob(\theta)-2\Lob(\pi/2-\theta)
\end{equation*}
(see for example \cite[Lemma 1]{Milnor:BULAM382}),
we have $4\Lob(\pi/6)=6\Lob(\pi/3)$, which equals the volume of
$S^3\setminus{E}$.
\end{proof}
\begin{proof}[Proof of Corollary~\ref{cor:Vol_CS}]
We see that
\begin{equation*}
  \Phi(u)
  =
  4H(u)-4\pi\sqrt{-1}u-4H(0)
\end{equation*}
and so
\begin{equation*}
  f(u)
  =
  H(u)-\pi\sqrt{-1}u-uv/4-\sqrt{-1}\Vol\left(S^3\setminus{E}\right)
\end{equation*}
since $\CS\left(S^3\setminus{E}\right)=0$.
Therefore from Theorem~\ref{thm:Yoshida}
\begin{equation*}
\begin{split}
  \Vol\left(E_{u}\right)+\sqrt{-1}\CS\left(E_{u}\right)
  &=
  \frac{H(u)}{\sqrt{-1}}
  -\pi{u}
  -\frac{uv}{4\sqrt{-1}}
  -\frac{\pi}{2}\lambda(\gamma).
\end{split}
\end{equation*}
\end{proof}
\section*{Appendix. Optimistic limit}
In this appendix, we will consider integral surgeries along the figure-eight
knot.
\par
Let $E_{p,1}$ be the $p$-surgery along the figure-eight knot with $p\in\Z$.
Then our parameter $u$ and $v$ satisfies $pu+v=2\pi\sqrt{-1}$.
From Corollary~\ref{cor:Vol_CS} we have
\begin{equation*}
  \Vol(E_{p,1})+\sqrt{-1}\CS(E_{p,1})
  \equiv
  \frac{H(u)}{\sqrt{-1}}-\pi{u}-\frac{uv}{4\sqrt{-1}}
  -\frac{\pi}{2}\lambda(\gamma)
  \mod{\pi^2\sqrt{-1}\Z}.
\end{equation*}
By \cite[(4.6)]{Meyerhoff/Neumann:COMMH92}, we have
\begin{equation*}
  \lambda(\gamma)
  =\frac{2\pi\sqrt{-1}}{p}-\frac{v}{p}
  =u.
\end{equation*}
(Note that we follow orientation conventions of
\cite{Neumann/Zagier:TOPOL85,Yoshida:INVEM85}
and so the sign of $v$ is opposite to that in \cite{Meyerhoff/Neumann:COMMH92}.)
Therefore
\begin{equation*}
\begin{split}
  &\Vol(E_{p,1})+\sqrt{-1}\CS(E_{p,1})
  \\
  &\quad\equiv
  {}-\sqrt{-1}\Li_2\left(y^{-1}\exp(-u)\right)
  +\sqrt{-1}\Li_2\left(y\exp(-u)\right)
  -\sqrt{-1}u\left\{\log(-y)+\pi\sqrt{-1}\right\}
  \\
  &\quad\quad
  -\pi{u}-\frac{u}{4\sqrt{-1}}(2\pi\sqrt{-1}-pu)
  -\frac{\pi}{2}u
  \\
  &\quad\equiv
  \frac{1}{\sqrt{-1}}
  \left\{
    \Li_2\left(y^{-1}\exp(-u)\right)
    -\Li_2\left(y\exp(-u)\right)
    +u\log(-y)-u\pi\sqrt{-1}+\frac{pu^2}{4}
  \right\}
  \\
  &\quad\quad\quad\mod{\pi^2\sqrt{-1}\Z}
\end{split}
\end{equation*}
and so
\begin{multline*}
  -\CS(E_{p,1})+\sqrt{-1}\Vol(E_{p,1})
  \\
  \equiv
  \Li_2\left(y^{-1}\exp(-u)\right)
  -\Li_2\left(y\exp(-u)\right)
  +u\log(-y)-u\pi\sqrt{-1}+\frac{pu^2}{4}
  \mod{\pi^2\Z}.
\end{multline*}
\par
Now we define
\begin{equation*}
\begin{split}
  V_p(\xi,\eta)
  &:=
  \Li_2\left(\xi^{-1}\eta^{-1}\right)
  -\Li_2\left(\xi\eta^{-1}\right)
  +(\log(-\xi))(\log{\eta})-\pi\sqrt{-1}\log\eta+\frac{p}{4}(\log{\eta})^2
  \\
  &=
  H(\xi,\eta)+\frac{p}{4}(\log{\eta})^2-2\pi\sqrt{-1}\log\eta.
\end{split}
\end{equation*}
We will calculate
$\partial\,V_p(\xi,\eta)/\partial{\xi}\bigr|_{\xi=y,\eta=m^2}$ and
$\partial\,V_p(\xi,\eta)/\partial{\eta}\bigr|_{\xi=y,\eta=m^2}$.
First we have from \eqref{eq:y}
\begin{equation*}
  \left.
    \frac{\partial\,V_p(\xi,\eta)}{\partial{\xi}}
  \right|_{\xi=y,\eta=m^2}
  =
  \left.
    \frac{\partial\,H(\xi,\eta)}{\partial{\xi}}
  \right|_{\xi=y,\eta=m^2}
  =
  0.
\end{equation*}
Next we have
\begin{equation*}
  \left.
    \frac{\partial\,V_p(\xi,\eta)}{\partial\eta}
  \right|_{\xi=y,\eta=m^2}
  =
  \left.
    \frac{\partial\,H(\xi,\eta)}{\partial\eta}
  \right|_{\xi=y,\eta=m^2}
  +
  \frac{p\log\left(m^2\right)}{2m^2}
  -\frac{2\pi\sqrt{-1}}{m^2}.
\end{equation*}
On the other hand
\begin{equation*}
\begin{split}
  \frac{v}{2}
  &=
  \frac{d\,H(u,m^2)}{du}
  -\pi\sqrt{-1}
  \\
  &=
  \left.
    \frac{\partial\,H(\xi,\eta))}{d\xi}
  \right|_{\xi=y,\eta=m^2}
  +
  \left.
    \frac{\partial\,H(\xi,\eta)}{d\eta}
  \right|_{\xi=y,\eta=m^2}
  \times
  \frac{d\,m^2}{du}
  -
  \pi\sqrt{-1}
  \\
  &=
  \exp{u}\times
  \left.
    \frac{\partial\,H(\xi,\eta)}{d\eta}
  \right|_{\xi=y,\eta=m^2}
  -
  \pi\sqrt{-1}.
\end{split}
\end{equation*}
Therefore
\begin{equation*}
  \left.
    \frac{\partial\,V_p(\xi,\eta)}{\partial\eta}
  \right|_{\xi=y,\eta=m^2}
  =
  \frac{v/2+\pi\sqrt{-1}}{\exp{u}}
  +\frac{pu}{2\exp{u}}
  -\frac{2\pi\sqrt{-1}}{\exp{u}}
  =0.
\end{equation*}
\par
In \cite{Murakami:SURIK00} the first author calculated a `fake' limit of the
Witten--Reshetikhin--Turaev invariants of the $(p,1)$-Dehn surgery of the
figure-eight knot and obtained the following observation.
\begin{obs*}
Let $p$ be an integer between $-100$ and $100$, and
$\tau_p$ the `optimistic' limit of the Witten--Reshetikhin--Turaev
invariant of the $(p,1)$-Dehn surgery along the figure-eight knot.
Then there exists $(\xi_0,\eta_0)$ such that the following equalities hold
numerically (up to 8 digits):
\begin{enumerate}
  \item
  $\dfrac{\partial\,V_p}{\partial\xi}(\xi_0,\eta_0)
  =\dfrac{\partial\,V_p}{\partial\eta}(\xi_0,\eta_0)=0$, and
  \item
  $\tau_p=V_p(\xi_0,\eta_0)=\CS(E_{p,1})+\sqrt{-1}\Vol(E_{p,1})$.
\end{enumerate}
\end{obs*}
\par
Note that the observation above was confirmed by the second author
\cite{Yokota:GTM02}.
Note also that the sign of $\CS$ is reversed since the orientation convention here is different.
\bibliography{mrabbrev,hitoshi}
\bibliographystyle{hamsplain}
\end{document}